\documentclass{amsart}
\usepackage{tikz}
\usepackage{hyperref}
\usepackage[backend=bibtex]{biblatex}
\newcommand{\Dfn}[1]{\emph{#1}}%
\newcommand{\fS}{\mathfrak{S}}%
\newcommand{\fD}{\mathfrak{D}}%
\DeclareMathOperator{\lrmax}{LRMAX}%
\DeclareMathOperator{\ldes}{ldes}%
\DeclareMathOperator{\blocks}{bl}%
\DeclareMathOperator{\returns}{ret}%
\DeclareMathOperator{\rises}{RISES}%
\DeclareMathOperator{\ldr}{ldr}%
\title[An involution on Dyck paths]%
{An involution on Dyck paths that preserves the rise composition and
  interchanges the number of returns and the position of the first
  double fall}%
\author{Martin Rubey}%
\email{\url{Martin.Rubey@tuwien.ac.at}}%
\address{Fakult\"at f\"ur Mathematik und Geoinformation, TU Wien,
  Austria}%

\begin{document}
\maketitle
\begin{abstract}
  Motivated by a recent paper of Adin, Bagno and Roichman, we present
  an involution on Dyck paths that preserves the rise composition and
  interchanges the number of returns and the position of the first
  double fall.
\end{abstract}

\section{Introduction}
In a recent paper of Adin, Bagno and Roichman~\cite{arXiv:1611.06979}
a new symmetry property of $321$-avoiding permutations is used to
demonstrate the Schur positivity of the quasi-symmetric function
associated with a certain set of permutations.  More precisely, the
authors show using a recursively defined
bijection\footnote{\url{http://www.findstat.org/MapsDatabase/Mp00124}}
that
\[
\sum_{\pi\in\fS_n(321)} \mathbf x^{\lrmax(\pi)} q^{\blocks(\pi)}%
=\sum_{\pi\in\fS_n(321)} \mathbf x^{\lrmax(\pi)}
q^{n-\ldes(\pi^{-1})}\text{, where}
\]
\begin{itemize}
\item $\fS_n(321)$ is the set of $321$-avoiding permutations of
$\{1,\dots,n\}$, 
\item $\lrmax(\pi)$ is the \emph{set} of positions of the
  left-to-right maxima of $\pi$, and
  $\mathbf x^L = \prod_{i\in L} x_i$,
\item
  $\blocks(\pi)=\lvert\{i: \forall j\leq i:\pi(j)\leq i\}\rvert$ is
  the number of blocks of $\pi$ and
\item $\ldes(\pi)=\max\{0\}\cup\{i: \pi(i)>\pi(i+1)\}$ is the
  position of the last descent of $\pi$.
\end{itemize}
We exhibit an iterative involution that proves the joint symmetry of
the two statistics linked by Adin, Bagno and Roichman.

Let us first transform the statistics involved into statistics on
Dyck paths, using a bijection due to
Krattenthaler~\cite{MR1868978}\footnote{\url{http://www.findstat.org/MapsDatabase/Mp00119}}.
To do so, it is convenient to define a Dyck path as a path in
$\mathbb N^2$ consisting of $(1,0)$ and $(0,1)$ steps beginning at
the origin, ending at $(n,n)$, and never going below the diagonal
$y=x$.  We refer to $(1,0)$-steps also as \Dfn{east steps} or
\Dfn{falls}, and to $(0,1)$-steps also as \Dfn{north steps} or
\Dfn{rises}, and the number $n$ as the \Dfn{semilength} of the path.
Krattenthaler's bijection maps
\begin{itemize}
\item $\fS_n(321)$ to the set of Dyck paths $\fD_n$ of semilength
  $n$,
\item $\lrmax(\pi)$ to $\rises(D)$, the set of $x$-coordinates
  (plus~1) of the north steps of the Dyck path,
\item $\blocks(\pi)$ to $\returns(D)$, the number of \Dfn{returns} of
  the path, that is, the number of east steps that end on the main
  diagonal,
\item $\ldes(\pi^{-1})$ to $\ldr(D)$, the $y$-coordinate of the
  middle point of the last double rise, that is, the last pair of two
  consecutive north steps.  For the path without double rises, we set
  $\ldr(D)=0$.
\end{itemize}

The purpose of this note is to provide an involution $\Phi$ that
implies
\begin{equation}
  \label{eq:switch}
  \sum_{\pi\in\fD_n} \mathbf x^{\rises(\pi)} p^{\returns(D)}
  q^{n-\ldr(D)}%
  =\sum_{\pi\in\fD_n} \mathbf x^{\rises(\pi)} p^{n-\ldr(D)}
  q^{\returns(D)}.
\end{equation}

\section{An involution on Dyck paths}

In the spirit
of~\cite{MR3422679}\footnote{\url{http://www.findstat.org/MapsDatabase/Mp00118}},
let us first introduce an invertible map $\phi$ from Dyck paths that
do not end with a single east step and satisfy $\ldr(D) < n-1$ (that
is, $\ldr(D)$ is not maximal), such that
\begin{itemize}
\item $\rises(\phi(D)) = \rises(D)$
\item $\returns(\phi(D)) = \returns(D)+1$
\item $\ldr(\phi(D)) = \ldr(D)+1$
\end{itemize}
Iterating this map (or its inverse) we obtain a map $\Phi$ that
proves Equation~\eqref{eq:switch}.

To define $\phi$, decompose $D$ as the concatenation of three paths
$P$, $Q$ and $R$ where
\begin{itemize}
\item $R$ is maximal such that $R$ contains only single rises, begins
  with a rise and $Q$ ends with an east step, and
\item $Q\,R$ is a prime Dyck path, that is, its only return is the
  final east step.
\end{itemize}
Let $n_0\, Q_1\, n_1\,\dots\, n_k\, Q_k$ be the prime decomposition
of $Q$, that is, all the $Q_i$ are prime Dyck paths and the $n_i$ are
sequences of north steps.  Note that $n_0$ is non-empty, because $Q$
does not return to the diagonal.  We now distinguish three cases:
\begin{enumerate}
\item If $n_1\neq\emptyset$, define
  $\phi(D)=P\,Q_1\, n_1\, \dots\, n_k\, Q_k\, n_0\, R$.
\item If $n_1=\dots=n_k=\emptyset$, let $n$ be the sequence of north
  steps of the final rise in $Q$ except one, and let $Q'$ be $Q$ with
  $n_0$ and $n$ erased.  Then define $\phi(D)=P\, n\, Q'\, n_0\, R$.
\item Otherwise, let $n$ be the sequence of north steps of the final
  rise in $Q$ except one, and let $Q'$ be $Q$ with $n_0$ and $Q_1$
  erased.  Then define $\phi(D)=P\, Q_1\, n\, Q'\, n_0\, R$.
\end{enumerate}
To show that $\phi$ is invertible we provide an inverse for each of
these three cases.  Let $D$ be a Dyck path that does not end with a
single east step and has at least two returns.  Decompose $D$ as the
concatenation of three paths $P$, $Q$ and $R$ where
\begin{itemize}
\item $R$ is maximal such that $R$ contains only single rises and
  begins with a rise, and
\item $Q\,e\,n\,R$ is a Dyck path with precisely two returns, $e$ is
  a maximal sequence of east steps and $n$ is a maximal sequence of
  north steps.
\end{itemize}
Let $n_0\, Q_1\, n_1\,\dots\, Q_k\, n_k$ be the prime
decomposition of $Q$.
\begin{enumerate}
\item If $Q$ ends with a double rise, we have
  $\phi^{(-1)}(D)=P\, n\, Q\, e\, R$.
\item Otherwise, if $Q$ has no returns, we have
  $\phi^{(-1)}(D)=P\, n\, Q_1\, n_1\,\dots\, Q_k\, n_k\, n_0\, e\,
  R$.
\item Otherwise,
  $\phi^{(-1)}(D)=P\, n\, Q_1\, Q_2\,n_2\,\dots\, Q_k\, n_k\, n_1\,
  e\, R$.
\end{enumerate}

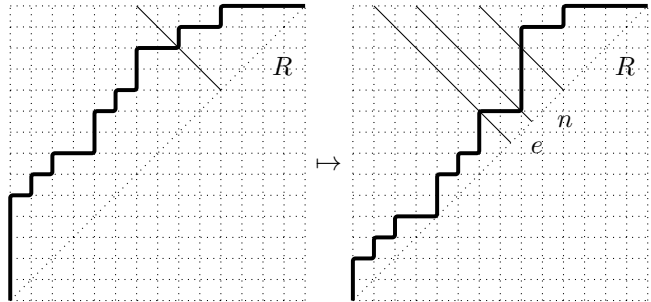
\begin{figure}
  \centering
  \begin{tikzpicture}[scale=0.28]
    \draw[dotted] (0, 0) grid (14, 14) -- (0,0);
    \draw[rounded corners=1, line width=1.5]%
    (0,0)--(0,5)--%
    (1,5)--(1,6)--%
    (2,6)--(2,7)--%
    (4,7)--(4,9)--%
    (5,9)--(5,10)--%
    (6,10)--(6,12)--%
    (8,12)--(8,13)--%
    (10,13)--(10,14)--%
    (14,14);
    \draw (6,14)--(10,10);
    \node[below right] at (12,12) {$R$};
  \end{tikzpicture}
  \raisebox{50pt}{$\mapsto$}
  \begin{tikzpicture}[scale=0.28]
    \draw[dotted] (0, 0) grid (14, 14) -- (0,0);
    \draw[rounded corners=1, line width=1.5]%
    (0,0)--(0,2)--%
    (1,2)--(1,3)--%
    (2,3)--(2,4)--%
    (4,4)--(4,6)--%
    (5,6)--(5,7)--%
    (6,7)--(6,9)--%
    (8,9)--(8,13)--%
    (10,13)--(10,14)--%
    (14,14);
    \draw (6,14)--(10,10);
    \draw (3,14)--(8.5,8.5);
    \draw (1,14)--(7.5,7.5);
    \node[below right] at (9.25,9.25) {$n$};
    \node[below right] at (8,8) {$e$};
    \node[below right] at (12,12) {$R$};
  \end{tikzpicture}
  \caption{A path and its image for case (1)}
\end{figure}

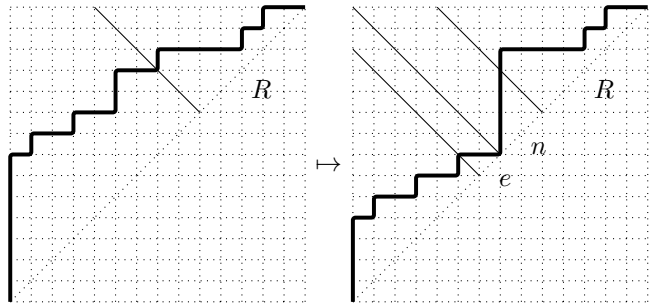
\begin{figure}
  \centering
  \begin{tikzpicture}[scale=0.28]
    \draw[dotted] (0, 0) grid (14, 14) -- (0,0);
    \draw[rounded corners=1, line width=1.5]%
    (0,0)--(0,7)--%
    (1,7)--(1,8)--%
    (3,8)--(3,9)--%
    (5,9)--(5,11)--%
    (7,11)--(7,12)--%
    (11,12)--(11,13)--%
    (12,13)--(12,14)--%
    (14,14);
    \draw (4,14)--(9,9);
    \node[below right] at (11,11) {$R$};
  \end{tikzpicture}
  \raisebox{50pt}{$\mapsto$}
  \begin{tikzpicture}[scale=0.28]
    \draw[dotted] (0, 0) grid (14, 14) -- (0,0);
    \draw[rounded corners=1, line width=1.5]%
    (0,0)--(0,4)--%
    (1,4)--(1,5)--%
    (3,5)--(3,6)--%
    (5,6)--(5,7)--%
    (7,7)--(7,12)--%
    (11,12)--(11,13)--%
    (12,13)--(12,14)--%
    (14,14);
    \draw (4,14)--(9,9);
    \draw (0,14)--(7,7);
    \draw (0,12)--(6,6);
    \node[below right] at (11,11) {$R$};
    \node[below right] at (8,8) {$n$};
    \node[below right] at (6.5,6.5) {$e$};
  \end{tikzpicture}
  \caption{A path and its image for case (2)}
\end{figure}

\begin{figure}
  \centering
  \begin{tikzpicture}[scale=0.28]
    \draw[dotted] (0, 0) grid (14, 14) -- (0,0);
    \draw[rounded corners=1, line width=1.5]%
    (0,0)--(0,5)--%
    (1,5)--(1,6)--%
    (2,6)--(2,7)--%
    (4,7)--(4,8)--%
    (5,8)--(5,11)--%
    (7,11)--(7,12)--%
    (8,12)--(8,13)--%
    (9,13)--(9,14)--%
    (14,14);
    \draw (4,14)--(9,9);
    \node[below right] at (11,11) {$R$};
  \end{tikzpicture}
  \raisebox{50pt}{$\mapsto$}
  \begin{tikzpicture}[scale=0.28]
    \draw[dotted] (0, 0) grid (14, 14) -- (0,0);
    \draw[rounded corners=1, line width=1.5]%
    (0,0)--(0,2)--%
    (1,2)--(1,3)--%
    (2,3)--(2,4)--%
    (4,4)--(4,7)--%
    (5,7)--(5,8)--%
    (7,8)--(7,12)--%
    (8,12)--(8,13)--%
    (9,13)--(9,14)--%
    (14,14);
    \draw (4,14)--(9,9);
    \draw (1,14)--(7.5,7.5);
    \draw (0,13)--(6.5,6.5);
    \node[below right] at (11,11) {$R$};
    \node[below right] at (8.5,8.5) {$n$};
    \node[below right] at (7,7) {$e$};
  \end{tikzpicture}
  \caption{A path and its image for case (3)}
\end{figure}
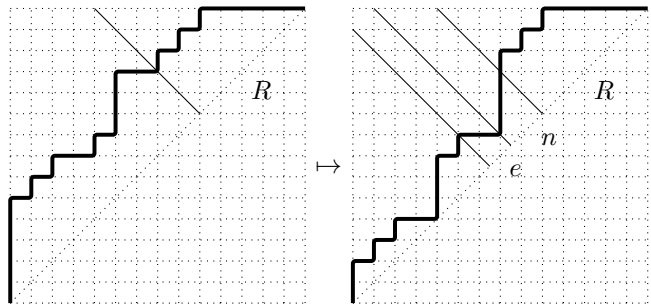

\printbibliography
\end{document}